\documentclass[10pt]{amsart}
\usepackage{epsfig}
\usepackage{amsmath, amssymb, euscript,  enumerate}
\usepackage{bbm}
\usepackage{color}
\usepackage{verbatim}

\newcommand{\dist}{{\text{\rm dist}}}

\newcommand{\supp}{\text{\rm supp}}

\newcommand{\ap}{\alpha}             
\newcommand{\bt}{\beta}
\newcommand{\gm}{\gamma}             \newcommand{\Gm}{\Gamma}
\newcommand{\dt}{\delta}             \newcommand{\Dt}{\Delta}
\newcommand{\vep}{\varepsilon}

\newcommand{\ld}{\lambda}            \newcommand{\Ld}{\Lambda}
             
\newcommand{\vp}{\varphi}
             \newcommand{\Om}{\Omega}
            \newcommand{\iy}{\infty}
            
\newcommand{\f}{\frac}

\newcommand{\fD}{{\mathfrak D}}

\newcommand{\fF}{{\mathfrak F}}

\newcommand{\fL}{{\mathfrak L}}

\newcommand{\fa}{{\mathfrak a}}

\newcommand{\fe}{{\mathfrak e}}

\newcommand{\BN}{{\mathbb N}}

\newcommand{\BR}{{\mathbb R}}

\newcommand{\cD}{{\mathcal D}}
\newcommand{\cE}{{\mathcal E}}
\newcommand{\cF}{{\mathcal F}}

\newcommand{\cI}{{\mathcal I}}
\newcommand{\cJ}{{\mathcal J}}
\newcommand{\cK}{{\mathcal K}}

\newcommand{\cQ}{{\mathcal Q}}

\newcommand{\cS}{{\mathcal S}}

\newcommand{\la}{\langle}          \newcommand{\ra}{\rangle}
\newcommand{\s}{\setminus}         
\newcommand{\n}{\nabla}            \newcommand{\e}{\eta}
\newcommand{\pa}{\partial}        
       
    \newcommand{\ds}{\displaystyle}

 \newcommand{\pf }{\noindent{\it Proof. }}
\newcommand{\rk }{\noindent{\it Remark. }}

\newcommand{\aee }{\text{\rm a.e.}} 
  \newcommand{\pv }{\text{\rm p.v.}}

\newcommand{\rL }{{\text{\rm L}}}  
   
\newcommand{\rH}{{\text{\rm H}}}
  
\newcommand{\rX}{{\text{\rm X}}}  \newcommand{\rY}{{\text{\rm Y}}}

\newcommand{\loc}{{\text{\rm loc}}}

\newtheorem{thm}[subsection]{Theorem}
\newtheorem{lemma}[subsection]{Lemma}
\newtheorem{cor}[subsection]{Corollary}

\newtheorem{prop}[subsection]{Proposition}

\newtheorem{defn}[subsection]{Definition}

\numberwithin{equation}{section}

\title[ Nonlocal Schr\"odinger operators with certain potential]{The Malgrange-Ehrenpreis theorem for nonlocal Schr\"odinger operators with certain potentials}
\author{ Woocheol Choi and Yong-Cheol Kim  }
\begin{document}
\begin{abstract} In this paper, we prove the Malgrange-Ehrenpreis theorem for nonlocal Schr\"odinger operators $L_K+V$ with nonnegative potentials $V\in L^q_{\loc}(\BR^n)$ for $q>\f{n}{2s}$ with $0<s<1$ and $n\ge 2$; that is to say, we obtain the existence of a fundamental solution $\fe_V$ for $L_K+V$ satisfying
\begin{equation*}\bigl(L_K+V\bigr)\fe_V=\dt_0\,\,\text{ in $\BR^n$ }\end{equation*}
in the distribution sense, where $\dt_0$ denotes the Dirac delta mass at the origin.
In addition, we obtain a decay of the fundamental solution $\fe_V$. 
\end{abstract}
\thanks {2000 Mathematics Subject Classification: 47G20, 45K05,
35J60, 35B65, 35D10 (60J75)}


\address{$\bullet$ Yong-Cheol Kim : Department of Mathematics Education, Korea University, Seoul 136-701,
Korea   }

\email{ychkim@korea.ac.kr}

\address{$\bullet$ Woocheol Choi : Department of Mathematics, Korea Institute for Advanced Study, Seoul 130-722, Korea } 

\email{wchoi@kias.re.kr}

\maketitle

\tableofcontents

\section{Introduction}
In 1954-1956, Malgrange and Ehrenpreis \cite{E1,E2,M} proved independently that any partial differential operator with constant coefficients which is not identically vanishing has  a fundamental solution in the space $\cD'(\BR^n)$ of distributions.
In this paper, we obtain an extention of their result to nonlocal Schr\"odinger operators with nonnegative potentials $V\in L^q_{\loc}(\BR^n)$ for $q>\f{n}{2s}$ with $0<s<1$ and $n\ge 2$.

We introduce integro-differential operators of the form
\begin{equation}L_K u(x)=\f{1}{2}\,\pv\int_{\BR^n}\mu(u,x,y)K(y)\,dy, \,\,x\in\BR^n,
\end{equation}
where $\mu(u,x,y)=2\,u(x)-u(x+y)-u(x-y)$ and the kernel $K:\BR^n\s\{0\}\to\BR_+$ satisfy the property
\begin{equation}\f{\ld\,c_{n,s}}{|y|^{n+2s}}\le K(y)=K(-y)\le\f{\Ld\,c_{n,s}}{|y|^{n+2s}},
\,\,s\in(0,1),\,0<\ld<\Ld<\iy,
\end{equation}
where $c_{n,s}$ is the normalization constant comparable to $s(1-s)$
given by
\begin{equation}c_{n,s}\int_{\BR^n}\f{1-\cos(\xi_1)}{|\xi|^{n+2s}}\,d\xi=1.
\end{equation}
Set $\fL=\{L_K:K\in\cK\}$ where $\cK$ denotes the family of all kernels $K$ satisfying (1.2). In particular, if $K(y)=c_{n,s}|y|^{-n-2s}$, then $L_K=(-\Dt)^s$ is the fractional Laplacian and it is well-known that
$$\lim_{s\to 1^-}(-\Dt)^s=-\Dt u$$ for any function $u$ in the Schwartz space $\cS(\BR^n)$.

In what follows, we consider the nonlocal Schr\"odinger operators given by
\begin{equation}\rL_V:=L_K+V\end{equation}
where $K\in\cK$ and $V\in L^q_{\loc}(\BR^n),$ $q>\f{n}{2s},$ $s\in(0,1),$ $n\ge 2,$ is a nonnegative potential. 
We are interested in the existence of a fundamental solution for the operator $\rL_V$. Let $\cD'(\BR^n)$ be the space of all distributions on $\BR^n$. Given $f\in\cD'(\BR^n)$, we say that a real-valued Lebesgue measurable function $u$ on $\BR^n$ satisfies the equation $\rL_V u=f$ {\it in the sense of $\cD'(\BR^n)$}, if the linear map $\rL_V u:C^{\iy}_c(\BR^n)\to\BR$  given by
\begin{equation*}\la\rL_V u,\vp\ra=\int_{\BR^n}u(y)\rL_V\vp(y)\,dy
\end{equation*}is well-defined and a distribution, and also
$\la\rL_V u,\vp\ra=\la f,\vp\ra$ for all $\vp\in C^{\iy}_c(\BR^n).$

\begin{thm} There exists a fundamental solution $\fe_V\in\cD'(\BR^n)$ for the nonlocal Schr\"odinger operator $\rL_V$, i.e. it satisfies that
$$\rL_V\fe_V=\dt_0\,\,\text{ in the sense of $\cD'(\BR^n)$, }$$
where $\dt_0$ is the Dirac delta mass at the origin. Moreover, there exists a universal constant $C>0$ depending on $n,s,\ld$ and $\Ld$ such that
\begin{equation}0\le\fe_V(x)\le\f{C}{|x|^{n-2s}}\,\,\text{ for any $x\in\BR^n\s\{0\}.$ }
\end{equation}
\end{thm}

\rk In case that $L_K=(-\Dt)^s$ and $V=0$, its fundamental solution was obtained in \cite{CS}.

The paper is organized as follows. In Section 2, we define several function spaces and give the fractional Sobolev embedding theorem which was proved in \cite{SV,DPV}.
In Section 3, we define weak solutions of the nonlocal equation $\rL_V=f$ in a bounded domain $\Om\subset\BR^n$ with Lipschitz boundary and obtain a relation between weak solutions (weak subsolutions, weak supersolutions) and minimizers (subminimizers, superminimizers) of the energy functional for the operator $\rL_V$, respectively.
In Section 4, we obtain a Rellich-Kondrachov compactness theorem, a weak maximum principle and a comparison principle for $\rL_V$.
In Section 5, we obtain an extension of the Malgrange-Ehrenpreis theorem to nonlocal Schr\"odinger operators $\rL_V$, that is, we furnish the proof of the existence of a fundamental solution $\fe_V$ for $\rL_V$ by using functional analysis stuffs. Moreover we obtain a decay of the fundamental solution $\fe_V$.

\section{Preliminaries}
Let $\cF^n$ be the family of all real-valued Lebesgue measurable functions on $\BR^n$.
Let $\Om$ be a bounded open domain in $\BR^n$ with {\it Lipschitz boundary}
and let $K\in\cK_s$. Let $\rX^s(\Om)$ be the function space of all $u\in\cF^n$ on $\BR^n$ such that $u|_\Om\in L^2(\Om)$
and
\begin{equation*}\iint_{\BR^{2n}_\Om}\f{|u(x)-u(y)|^2}{|x-y|^{n+2s}}\,dx\,dy<\iy
\end{equation*}
where $\BR^{2n}_D:=(\BR^n\times\BR^n)\s(D^c\times D^c)$ for a set $D\subset\BR^n$. We
also denote by
\begin{equation}\rX^s_0(\Om)=\{v\in\rX^s(\Om):u=0\,\,\aee\text{ in $\BR^n\s\Om$ }\}
\end{equation}
Note that $\rX^s(\Om)$ and $\rX^s_0(\Om)$ are not
empty, because $C^2_0(\Om)\subset\rX^s_0(\Om)$. Then we see that $(\rX^s(\Om),\|\cdot\|_{\rX^s(\Om)})$ is a normed space,
where the norm $\|\cdot\|_{\rX^s(\Om)}$ is defined by
\begin{equation}\|u\|_{\rX^s(\Om)}:=\|u\|_{L^2(\Om)}+\biggl(\iint_{\BR^{2n}_\Om}\f{|u(x)-u(y)|^2}{|x-y|^{n+2s}}\,dx\,dy\biggr)^{1/2}<\iy.
\end{equation}
For $p\ge 1$ and $s\in(0,1)$, let $W^{s,p}(\Om)$ be the usual fractional Sobolev spaces
with the norm
\begin{equation}\|u\|_{W^{s,p}(\Om)}:=\|u\|_{L^p(\Om)}+[u]_{W^{s,p}(\Om)}<\iy
\end{equation}
where the seminorm $[\,\cdot\,]_{W^{s,p}(\Om)}$ is defined by
$$[u]_{W^{s,p}(\Om)}=\biggl(\iint_{\Om\times\Om}
\f{|u(x)-u(y)|^p}{|x-y|^{n+sp}}\,dx\,dy\biggr)^{1/p}.$$ In what follows, we write $H^s(\Om)=W^{s,2}(\Om)$.
When $\Om=\BR^n$ in (2.3), we can similarly define the spaces $W^{s,p}(\BR^n)$  and $H^s(\BR^n)=W^{s,2}(\BR^n)$ for $s\in(0,1)$. 

By \cite{SV}, there exists a constant $c>1$ depending only on $n,s$ and $\Om$ such that
\begin{equation}\|u\|_{\rX^s_0(\Om)}\le\|u\|_{\rX^s(\Om)}\le c\,\|u\|_{\rX^s_0(\Om)}
\end{equation}
for any $u\in\rX^s_0(\Om)$, where
\begin{equation}\|u\|_{\rX^s_0(\Om)}:=\biggl(\iint_{\BR^{2n}_\Om}\f{|u(x)-u(y)|^2}{|x-y|^{n+2s}}\,dx\,dy\biggr)^{1/2}.
\end{equation} Thus $\|\cdot\|_{\rX^s_0(\Om)}$ is a norm on $\rX^s_0(\Om)$ which is equivalent to (2.2). Moreover it is known \cite{SV} that $(\rX^s_0(\Om),\|\cdot\|_{\rX^s_0(\Om)})$ is a
Hilbert space with inner product
\begin{equation}\la
u,v\ra_{\rX^s_0(\Om)}:=\iint_{\BR^{2n}_\Om}\f{(u(x)-u(y))(v(x)-v(y))}{|x-y|^{n+2s}}\,dx\,dy.
\end{equation}
Let $\rX^s_0(\Om)^*$ be the dual space of $\rX^s_0(\Om)$; that is, the family of all
bounded linear functionals on $\rX^s_0(\Om)$. Then we know that
$(\rX^s_0(\Om)^*,\|\cdot\|_{\rX^s_0(\Om)^*})$ is a Hilbert space, where the norm
$\|\cdot\|_{\rX^s_0(\Om)^*}$ is given by
$$\|u\|_{\rX^s_0(\Om)^*}:=\sup\{u(v):v\in\rX^s_0(\Om), \|v\|_{\rX^s_0(\Om)}\le
1\},\,\,u\in\rX^s_0(\Om)^*.$$ 
For $s>0$, we consider the class $\rH^s(\BR^n)=\{u\in L^2(\BR^n):|\xi|^{2s}|\widehat u(\xi)|\text{ is in $L^2(\BR^n)$ }\}$ whose norm is defined by 
$$\|u\|_{\rH^s(\BR^n)}=\biggl(\int_{\BR^n}(1+|\xi|^s)^2|\widehat u(\xi)|^2\,d\xi\biggr)^{1/2},$$ 
where the Fourier transform of $u$ is defined by
$$\widehat u(\xi)=\int_{\BR^n}e^{-i\la x,\xi\ra}u(x)\,dx.$$
Then, by the Plancherel theorem, it is easy to check that $\rH^s(\BR^n)=H^s(\BR^n)$ and they are norm-equivalent. 

We also define the {\it homogeneous fractional Sobolev spaces} $\dot H^s(\BR^n)$ 
by the closure of $\cS(\BR^n)$ with respect to the norm
\begin{equation*}\|u\|_{\dot H^s(\BR^n)}=\biggl(\int_{\BR^n}|\xi|^{2s}|\widehat u(\xi)|^2\,d\xi\biggr)^{1/2}.
\end{equation*}
From a direct calculation \cite{DPV}, it turns out that
\begin{equation}\|u\|_{\dot H^s(\BR^n)}=2\,c^{-1}_{n,s}\iint_{\BR^n\times\BR^n}\f{|u(x)-u(y)|^2}{|x-y|^{n+2s}}\,dx\,dy
\end{equation}  where $c_{n,s}$ is the universal constant given in (1.3).

\begin{lemma} $(a)$ If $f\in\rX^s_0(\Om)$, then $f\in H^s(\BR^n)$ and moreover 
$$\f{c_{n,s}}{2}\,\|f\|_{\dot H^s(\BR^n)}=\|f\|_{\rX^s_0(\Om)}\le\|f\|_{H^s(\BR^n)}=\|f\|_{\rX^s(\Om)}\le c\,\|f\|_{\rX^s_0(\Om)} $$
where $c>1$ is the constant given in $(2.4)$.

$(b)$ If $f\in H^s(\BR^n)$, then $\|f\|_{\rX^s_0(\Om)}\le\|f\|_{\rX^s(\Om)}\le\|f\|_{H^s(\BR^n)}$.
\end{lemma}

\pf (a) It easily follows from (2.7) and Lemma 5 in \cite{SV}. (b) It is also very straightforward. \qed

\,\,\,Next we state the fractional Sobolev embedding theorem, which was proved in \cite{SV,DPV}.

\begin{prop} If $0\le s<\f{n}{2}$, then the space $\dot H^s(\BR^n)$ is continuously embedded in $L^{\f{2n}{n-2s}}(\BR^n)$, i.e. there is a constant $C>0$ depending only on  $n,s$ such that
\begin{equation*}\|u\|_{L^{\f{2n}{n-2s}}(\BR^n)}\le C\,\|u\|_{\dot H^s(\BR^n)}.
\end{equation*}
\end{prop}

\section{Weak solutions and minimizers for $L_K+V$ }

In this section, we define weak solutions, weak subsolutions and supersolutions of the nonlocal equation $L_K u+V u=0$ on $\Om$. To comprehend them well, we obtain a relation between weak solutions (weak subsolutions and weak supersolutions) and minimizers (subminimizers and superminimizers) of the energy functional for the nonlocal operator $L_K+V$, respectively.

From now on, we always assume that {\it $V$ is a nonnegative potential in $L^q_{\loc}(\BR^n)$   for $q>\f{n}{2s}$ with $s\in(0,1)$ and $n\ge 2$.}
We denote by $L^2_V(\Om)$ the weighted $L^2$ class of all real-valued measurable functions $g$ on $\BR^n$ satisfying
$$\|g\|^2_{L^2_V(\Om)}:=\int_{\Om}|g(y)|^2\,V(y)\,dy<\iy.$$
We consider a bilinear form defined by
\begin{equation*}\la u,v\ra_K=\iint_{\BR^n\times\BR^n}(u(x)-u(y))(v(x)-v(y))\,K(x-y)\,dx\,dy\,\,\text{ for $u,v\in\rX^s(\Om)$.}
\end{equation*}
From (b) of Corollary 5.2 below, we see in advance that
\begin{equation}\la u,v\ra_K=\int_{\BR^n}L_K u(y)\,v(y)\,dy=\int_{\BR^n}u(y)\,L_K v(y)\,dy\sim\la u,v\ra_{\rX^s_0(\Om)}
\end{equation} for $u,v\in\rX^s_0(\Om)$.
\begin{defn} Let $V\in L^q_{\loc}(\BR^n)$ for $q>\f{n}{2s}$ with $s\in(0,1)$ and $n\ge 2$. Then we say that a function $u\in\rX^s_0(\Om)$ is  a {\rm weak solution} of the nonlocal equation 
\begin{equation}\rL_V u:=L_K u+V u=f\,\,\text{ in $\Om$}
\end{equation} where $f\in\rX^s_0(\Om)^*$,
if it satisfies the weak formulation
\begin{equation}\int_{\BR^n}L_K u(y)\,\vp(y)\,dy+\int_{\BR^n}V(y)u(y)\vp(y)\,dy=\la f,\vp\ra\end{equation}
for any $\vp\in\rX^s_0(\Om)$. Here $\la\,\cdot\,,\,\cdot\,\ra$ denotes the pair between $\rX^s_0(\Om)^*$ and $\rX^s_0(\Om)$.
\end{defn}

\,In fact, it turns out that the weak solution of the equation (3.2) is the minimizer of the energy functional
\begin{equation}\cE_V(v)=\|v\|^2_{\rX^s_0(\Om)}+\|v\|^2_{L^2_V(\Om)}-2\la f,v\ra,\,\,v\in\rY^s_0(\Om):=\rX^s_0(\Om)\cap L^2_V(\Om),
\end{equation}
where $\rY^s_0(\Om)$ be a Hilbert subspace of $\rX^s_0(\Om)$ which is endowed with the norm
$$\|u\|_{\rY^s_0(\Om)}:=\sqrt{\|u\|^2_{\rX^s_0(\Om)}+\|u\|^2_{L^2_V(\Om)}}<\iy,\,\,u\in\rY^s_0(\Om).$$
We consider function spaces $\rY^s_0(\Om)^+$ and $\rY^s_0(\Om)^-$ defined by $$\rY^s_0(\Om)^\pm=\{v\in\rY^s(\Om):v^{\pm}\in\rX^s_0(\Om)\}$$ 
where $\rY^s(\Om):=\rX^s(\Om)\cap L^2_V(\Om)$ be the subspace of $\rX^s(\Om)$ which is endowed with the norm
$$\|u\|_{\rY^s(\Om)}:=\sqrt{\|u\|^2_{\rX^s(\Om)}+\|u\|^2_{L^2_V(\Om)}}<\iy,\,\,u\in\rY^s(\Om).$$
Then we see that $\rX^s_0(\Om)=\rY^s_0(\Om)^+\cap\rY^s_0(\Om)^-$. 

\,We now define weak subsolutions and supersolutions of the nonlocal equation (3.2) as follows. 

\begin{defn} We say a function $u\in\rY^s_0(\Om)^-\,(\,\rY^s_0(\Om)^+\,)$ is a {\rm weak subsolution (weak supersolution)} of the nonlocal equation $(3.2)$, if it satisfies that
\begin{equation}\int_{\BR^n}L_K u(y)\,\vp(y)\,dy+\int_{\Om}V(y)u(y)\vp(y)\,dy\le (\,\ge\,)\,\la f,\vp\ra
\end{equation} for every nonnegative $\vp\in\rX^s_0(\Om)$. Also we say that a function $u$ is a {\rm weak solution} of the nonlocal equation $(3.2)$, if it is both a weak subsolution and a weak supersolution. So any weak solution $u$ of the equation $(3.2)$ must be in $\rX^s_0(\Om)$ and satisfies $(3.3)$.
\end{defn}

In the next, we furnish the definition of subminimizer and superminimizer of the functional (3.4) to get better understanding of weak subsolutions and supersolutions of the nonlocal equation (3.2).

\begin{defn} $(a)$ We say that a function $u\in\rY^s_0(\Om)^-$ is a {\rm subminimizer} of the functional $(3.4)$ over $\rY^s_0(\Om)^-$, if it satisfies that
\begin{equation}\cE_V(u)\le\cE_V(u+\vp)
\end{equation} for all nonpositive $\vp\in\rX^s_0(\Om)$. Also we say that a function $u\in\rY^s_0(\Om)^+$ is a {\rm superminimizer} of the functional $(3.4)$ over $\rY^s_0(\Om)^+$, if it satisfies $(3.6)$ for all nonnegative $\vp\in\rX^s_0(\Om)$. 

$(b)$ We say that a function $u$ is a {\rm minimizer} of the functional $(3.4)$ over $\rX^s_0(\Om)$, if it is both a subminimizer and a superminimizer. So any minimizer $u$ must be in $\rX^s_0(\Om)$ and satisfies $(3.6)$ for all $\vp\in\rX^s_0(\Om)$.
\end{defn}

\begin{lemma} If $s\in(0,1)$, then there is a unique minimizer of the functional $(3.4)$. Moreover, a function $u\in\rY^s_0(\Om)^-\,(\,\rY^s_0(\Om)^+\,)$ is a subminimizer $(superminimizer)$ of the functional $(3.4)$ over $\rY^s_0(\Om)^-\,(\,\rY^s_0(\Om)^+\,)$ if and only if it is a weak subsolution $(weak$ $supersolution)$ of the nonlocal equation $(3.2)$. In particular, a function $u\in\rX^s_0(\Om)$ is a minimizer of the functional $(3.4)$ if and only if it is a weak solution of the equation $(3.2)$.
\end{lemma}

\pf Using standard method of calculus of variation, we proceed with our proof. We now take any minimizing sequence $\{u_k\}\subset\rX^s_0(\Om)$. By applying Theorem 4.1 below, we can take a subsequence $\{u_{k_j}\}\subset\rX^s_0(\Om)$ converging strongly to $u\in L^2(\Om)$. So there exists a subsequence $\{u_{k_i}\}$ of $\{u_{k_j}\}$ which converges $\aee$ in $\Om$ to $u\in\rX^s_0(\Om)$. Thus by applying Fatou's lemma we can show that the energy functional $\cE_V$ is weakly semicontinuous in $\rX^s_0(\Om)$. This implies that $u$ is a minimizer of (3.4). Its uniqueness also follows from the strict convexity of the functional (3.4). 

Next, we show the equivalency only for the weak supersolution case, because the other case can be done in a similar way. First, if $u\in\rY^s_0(\Om)^+$, then we observe that
\begin{equation}\begin{split}
\cE_V(u+\vp)-\cE_V(u)&=2\la u,\vp\ra_K+\int_{\Om}V(y)u(y)\vp(y)\,dy-2\la f,\vp\ra\\
&\,\qquad\qquad\qquad+\|\vp\|^2_{\rX^s_0(\Om)}+\|\vp\|^2_{L^2_V(\Om)}
\end{split}\end{equation} for all nonnegative $\vp\in\rX^s_0((\Om)$. This implies that a weak supersolution $u\in\rY_0^s(\Om)^+$ of the equation (3.2) is a superminimizer of the functional (3.4) over $\rY_0^s(\Om)^+$. 

On the other hand, we suppose that $u\in\rY^s_0(\Om)^+$ is a superminimizer 
of the functional (3.4). Then it follows from (3.7) that
$$2\la u,\vp\ra_K+2\int_{\Om}V(y)u(y)\vp(y)\,dy-2\la f,\vp\ra+\|\vp\|^2_{\rX^s_0(\Om)}+\|\vp\|^2_{L^2_V(\Om)}\ge 0$$ for all nonnegative $\vp\in\rX^s_0(\Om)$. Since $\vep\vp\in\rX^s_0(\Om)$ and it is nonnegative for any $\vep>0$ and $\vp\in\rX^s_0(\Om)$, we obtain that
$$2\la u,\vp\ra_K+2\int_{\Om}V(y)u(y)\vp(y)\,dy-2\la f,\vp\ra+\vep\|\vp\|^2_{\rX^s_0(\Om)}+\vep\|\vp\|^2_{L^2_V(\Om)}\ge 0$$ for any $\vep>0$. Taking $\vep\to\iy$, we can conclude that
$$\la u,\vp\ra_K+\int_{\Om}V(y)u(y)\vp(y)\,dy-2\la f,\vp\ra\ge 0$$ for any nonnegative $\vp\in\rX^s_0(\Om)$.
Hence $u$ is a weak supersolution of the equation (3.2). Therefore we are done. \qed

\section{Rellich-Kondrachov Compactness Theorem, a Weak Maximum Principle and a Comparison Principle for $L_K+V$ }

In this section, we obtain the {\it Rellich-Kondrachov compactness theorem}, a {\it weak maximum principle} and a {\it comparison principle} for $L_K+V$ which will play a crucial role in proving the existence of a fundamental solution for the nonlocal Schr\"odinger operators in the next section.

First we get a compactness result $\rY^s_0(\Om)\hookrightarrow L^2(\Om)$ by using the fact that  $\rX^s_0(\Om)$ and $\rY^s_0(\Om)$ are norm-equivalent and the precompactness of $\rY^s_0(\Om)$ in $L^2(\Om)$.

\begin{thm} Let $n\ge 1$, $s\in(0,1)$ and $2s<n$. If $\,u\in\rY^s_0(\Om)$, then there is a universal constant $C>0$ depending on $n,s$ and $\Om$ such that
$$\|u\|_{L^2(\Om)}\le C\,\|u\|_{\rY^s_0(\Om)}.$$
Moreover, any bounded sequence in $\rY^s_0(\Om)$ is precompact in $L^2(\Om)$.
\end{thm}

\pf We observe that $\rY^s_0(\Om)=\rX^s_0(\Om)$ and they are norm-equivalent, because it follows from Lemma 2.1 and the fractional Sobolev inequality \cite{DPV} that 
\begin{equation*}\|u\|_{L^2_V(\Om)}\le C\,\|V\|^{1/2}_{L^q(\Om)}\,\|u\|_{\rX^s_0(\Om)}
\end{equation*} with a universal constant $C>0$ depending on $n,s$ and $\Om$, and so 
\begin{equation}\|u\|^2_{\rX^s_0(\Om)}\le\|u\|^2_{\rY^s_0(\Om)}\le\bigl(1+C^2\|V\|_{L^q(\Om)}\bigr)\,\,\|u\|^2_{\rX^s_0(\Om)}
\end{equation} for any $u\in\rX^s_0(\Om)$. Applying Lemma 2.1 and the fractional Sobolev inequality again, we conclude that
\begin{equation*}\|u\|_{L^2(\Om)}\le C\,\|u\|_{H^s(\Om)}\le C\,\|u\|_{\rX^s_0(\Om)}\le C\,\|u\|_{\rY^s_0(\Om)}.
\end{equation*}

For the proof of the second part, take any bounded sequence $\{u_k\}$ in $\rY^s_0(\Om)$. Then it is also a bounded sequence in $\rX^s_0(\Om)$. Thus by Lemma 8 \cite{SV} there is a subsequence $\{u_{k_j}\}$ of the sequence and $u\in L^2(\Om)$ such that $u_{k_j}\to u$ in $L^2(\Om)$ as $j\to\iy$. Hence we complete the proof. \qed

\,\,\,Next, we give a weak maximum principle and a comparison principle for the nonlocal Schr\"odinger operators $\rL_V$
as follows.

\begin{thm} If $u$ is a weak supersolution of the nonlocal equation $\rL_V u=0$ in $\Om$ such that and $u\ge 0$ in $\BR^n\s\Om$, then $u\ge 0$ in $\Om$.
\end{thm}

\pf From the assumption, we see that $u^-=0$ in $\BR^n\s\Om$ and $u^+\in\rX^s_0(\Om)$, and thus $u^-,u\in\rX^s_0(\Om)$. 
Since $u^+\,u^-=0$ in $\BR^n$ and $u^+(x)u^-(y)\ge 0$ for all $x,y\in\BR^n$, we have that  
\begin{equation}\begin{split}
0&\le\la u,u^-\ra_{\rX^s_0(\Om)}+\int_{\Om}V(y)u(y)u^-(y)\,dy\\
&=\la u^+-u^-,u^-\ra_{\rX^s_0(\Om)}-\int_{\Om}V(y)[u^-(y)]^2\,dy\\
&\le-\|u^-\|^2_{\rX^s_0(\Om)}+\iint_{\BR^{2n}_{\Om}}\bigl(u^+(x)-u^+(y)\bigr)\bigl(u^-(x)-u^-(y)\bigr)K(x-y)\,dx\,dy\\
&=-\|u^-\|^2_{\rX^s_0(\Om)}-\iint_{\BR^{2n}_{\Om}}\bigl(u^+(x)u^-(y)+u^-(x)u^+(y)\bigr)K(x-y)\,dx\,dy\\
&\le-\|u^-\|^2_{\rX^s_0(\Om)}.
\end{split}\end{equation} This implies that $u^-=0$ in $\BR^n$, and hence $u\ge 0$ in $\Om$. \qed

\begin{cor} If $u$ is a weak subsolution of the nonlocal equation $\rL_V u=0$ in $\Om$ such that and $u\le 0$ in $\BR^n\s\Om$, then $u\le 0$ in $\Om$.
\end{cor}
 
\begin{cor} If $u$ is a weak subsolution and $v$ is a weak supersolution of the nonlocal equation $(3.2)$ such that $u\le v$ in $\BR^n\s\Om$, then $u\le v$ in $\Om$.
\end{cor}

\pf It immediately follows from Theorem 4.2. \qed

\section{Proof of the Malgrange-Ehrenpreis theorem for $L_K+V$ }

In this section, we study the existence of a fundamental solution $\fe_V$ for the nonlocal Schr\"odinger operators $\rL_V$, where $V$ is a nonnegative potential with $V\in L^q_{\loc}(\BR^n)$ for $q>\f{n}{2s}$ and $s\in(0,1)$ and $n\ge 2$. 

Let $T$ be an unbounded densely defined linear operator with domain $\fD(T)$ in a Hilbert space $H$ with the inner product $\la \cdot,\cdot\ra_H$. We denote by $\fD(T^*)$ the class of 
$\e\in H$ for which there exists a $\nu\in H$ such that
$$\la T(\upsilon),\e\ra_H=\la\upsilon,\nu\ra_H\,\,\text{ for all $\upsilon\in\fD(T)$.}$$
For each $\e\in\fD(T^*)$, we define $T^*(\e)=\nu$. Then we call $T^*$ the {\it adjoint} of $T$.

Let $\Gm(T)$ be the {\it graph} of such an operator $T$; that is, it is the linear subspace
$$\Gm(T)=\{(u,v)\in H\times H: u\in\fD(T)\text{ and } v=T(u)\}.$$
The operator $T$ is said to be {\it closed}, if $\Gm(T)$ is closed in $H\times H$. Also the operator $T$ is said to be {\it closable}, if there is a closed extension $T_0$ of $T$; that is, there is a closed operator $T_0$ with $T\subset T_0$ ( i.e. $\Gm(T)\subset\Gm(T_0)$ ). We call $\overline T$ the {\it closure} of $T$, i.e. the smallest closed extension of $T$.
Then the following two facts are well-known \cite{RS}:
\begin{equation}\begin{split}\,\,&(a)\,\,\text{ If $T$ is closable, then $\Gm(\overline T)=\overline{\Gm(T)}$ and $\overline T=T^{\ast\ast}$. }\\
\,\,&(b)\,\,\text{ If $T_1,T_2$ are densely defined operators with $T_1\subset T_2$, then $T_2^*\subset T_1^*$. }
\end{split}\end{equation}

For $s\in(0,1)$, we denote by $\rX^s_c:=\rX^s_c(\BR^n)$ the class of all $u\in\cF^n$ such that $u\in\rX^s_0(\Om)$ for some bounded domain $\Om$ in $\BR^n$. Similarly, we can define $\rY^s_c:=\rY^s_c(\BR^n)$.

\begin{lemma} If $K\in\cK_s$ for $s\in(0,1)$, then the operator $L_K:L^2(\BR^n)\to L^2(\BR^n)$ is a densely defined operator with domain $\fD(L_K)=C_c^{\iy}(\BR^n)$ and it is positive and symmetric. 
Moreover, there exists a unique closure $\overline{L_K}$ of $L_K$ which is self-adjoint and $\fD(\overline{L_K})=H^{2s}(\BR^n)$.
\end{lemma}

\pf Note that $C^{\iy}_c(\BR^n)$ is dense in $L^2(\BR^n)$ and $C^{\iy}_c(\BR^n)\subset\rX^s_c(\BR^n)$. By Theorem 4.1, it is easy to check that $L_K:L^2(\BR^n)\to L^2(\BR^n)$ is a densely defined operator with domain $\fD(L_K)=C^{\iy}_c(\BR^n)$. Also we see that $L_K$ is positive and symmetric, because $\la L_K u,u\ra_{L^2(\BR^n)}=\|u\|^2_{\rX^s_0(\BR^n)}\ge 0$ and $\la L_K u,v\ra_{L^2(\BR^n)}=\la u,L_K v\ra_{L^2(\BR^n)}$ for any $u,v\in C^{\iy}_c(\BR^n)$. 

To prove the existence of the closure $\overline {L_K}$ which is self-adjoint, it suffices to show that $[L_K\pm i](C^{\iy}_c(\BR^n))$ is dense in $L^2(\BR^n)$ by verifying that its orthogonal complement is $[L_K\pm i](C^{\iy}_c(\BR^n))^\perp=0$ (refer to p.257, \cite{RS}). Indeed, let us take any $\vp\in L^2(\BR^n)$ satisfying \begin{equation}\la\vp,[L_K\pm i]\phi\ra_{L^2(\BR^n)}=0\,\,\text{ for all $\phi\in\fD(L_K)=C_c^{\iy}(\BR^n)$. }
\end{equation}
Since the Fourier transform $\fF(f)=\widehat f$ is a homeomorphism from $L^2(\BR^n)$ onto itself, we have that
\begin{equation}\begin{split}0&=\la\vp,[L_K\pm i]\phi\ra_{L^2(\BR^n)}=\la\widehat\vp,\fF([L_K\pm i]\phi)\ra_{L^2(\BR^n)}\\&=\la\widehat\vp,(m(\xi)\pm i)\widehat\phi\,\ra_{L^2(\BR^n)}=\la(m(\xi)\mp i)\widehat\vp,\widehat\phi\,\ra_{L^2(\BR^n)}
\end{split}\end{equation}  for all $\phi\in\fD(L_K)=C_c^{\iy}(\BR^n)$, where $m(\xi)$ is the nonnegative function given by
\begin{equation}m(\xi)=\int_{\BR^n}(1-\cos\la y,\xi\ra)K(y)\,dy.
\end{equation}
Since $\fF(C_c^{\iy}(\BR^n))$ is dense in $L^2(\BR^n)$, it follows from the Plancherel theorem and (5.3) that $(m(\xi)\mp i)\widehat\vp=0$, and thus $\vp=0$.  

Next, we show the uniqueness of the closure $\overline{L_K}$ which is self-adjoint.  
Indeed, if $S$ is another self-adjoint closed extension of $L_K$, then we see that $\overline{L_K}\subset S$. Conversely, it follows from (5.1) that $[L_K]^{**}\subset S$, and hence $$S=S^*\subset[\overline{L_K}]^*=\overline{L_K}.$$

Finally, we consider the multiplication operator $M_0:L^2(\BR^n)\to L^2(\BR^n)$ with $\fD(M_0)=\fF(C_c^{\iy}(\BR^n))$ defined by $M_0(\widehat\vp)=m\,\widehat\vp$, where $m$ is the function in (5.4).
Then, by the Plancherel theorem, we see that $L_K$ is unitarily equivalent to $M_0$.
As in (5.3), we can also prove that there is a unique closure $\overline M_0$ of $M_0$ which is self-adjoint and $$\fD(\overline M_0)=\{\phi\in L^2(\BR^n):m\phi\in L^2(\BR^n)\}.$$ Thus we have that $\fD(\overline M_0)=H^{2s}(\BR^n)$, because we easily obtain that $$\ld|\xi|^{2s}\le m(\xi)\le\Ld|\xi|^{2s}$$ by (1.2) and (1.3). Since the Fourier transform $\cF$ is a unitary isomorphism from $L^2(\BR^n)$ to $L^2(\BR^n)$, the closure 
$\overline{L_K}$ of $L_K$ is unitarily equivalent to $\overline M_0$. Therefore we conclude that  $\fD(\overline{L_K})=H^{2s}(\BR^n)$. \qed

\,\,\,As by-products of Lemma 5.1, we get a very useful {\it nonlocal version of integration by parts} and the norm equivalence between the space $\rX^s_0(\Om)$ and the usual fractional Sobolev space $H^s(\BR^n)$ on a dense subspace $C_c^{\iy}(\Om)$ of $\rX^s_0(\Om)$.

\begin{cor} $(a)$  If $K\in\cK$ for $s\in(0,1)$, then there is a unique positive self-adjoint square root operator $Q$ of $\overline{L_K}$, i.e. $Q\circ Q=\overline{L_K}$. Also, it satisfies that
\begin{equation}\la L_K u,v\ra_{L^2(\BR^n)}=\la Q u,Q v\ra_{L^2(\BR^n)}\,\,\text{ and }\,\,\widehat{Q u}(\xi)=\sqrt{m(\xi)}\,\widehat u(\xi)
\end{equation} for any $u,v\in C^{\iy}_c(\BR^n)$, where $m(\xi)$ is the multiplier of $L_K$ given in $(5.4)$.

\,\,$(b)$ $\ds\int_{\Om}L_K u(y)\,v(y)\,dy=\int_{\Om}u(y)\,L_K v(y)\,dy=\la u,v\ra_{\rX^s_0(\Om)}$ for all $u,v\in\rX^s_0(\Om)$.

\,\,$(c)$ $\|u\|_{\rX^s_0(\Om)}\sim\|Q u\|_{L^2(\BR^n)}\sim\|u\|_{H^s(\BR^n)}$ for all $u\in C_c^{\iy}(\Om)$.
\end{cor}

\rk If $K(y)=c_{n,s}|y|^{-n-2s}$ for $s\in(0,1)$, then we note that $L_K=(-\Dt)^s$ and $Q=(-\Dt)^{s/2}$.

\,\pf (a) It can be obtained by Theorem 13.31 \cite{R}. (b) It can be shown by simple calculation. (c) It easily follows from (5.5) and (b). Therefore we complete the proof. \qed

\,\,\,Let $H_{V^2}(\Om)=L^2_{V^2}(\Om)\cap L^2(\Om)$ be a Hilbert space with inner product $$\la u,v\ra_{H_{V^2}(\Om)}=\la Vu,Vv\ra_{L^2(\Om)}+\la u,v\ra_{L^2(\Om)}.$$

\begin{lemma} $(a)$ The multiplication operator $M_V:L^2(\Om)\to L^2(\Om)$ with domain $\fD(M_V)=H_{V^2}(\Om)$ given by $M_V(u)=Vu$ is positive and self-adjoint. 

$(b)$ There is a unique positive self-adjoint square root operator $P$ of $M_V$, i.e.
$P\circ P=M_V$. 
\end{lemma}

\pf (a) For any fixed $v\in L^2(\Om)$, it is quite easy to check that the linear map $T_0:H_{V^2}(\Om)\to\BR$ defined by $T_0(u)=\la M_V u,v\ra_{L^2(\Om)}$ satifies that $$|T_0(u)|\le\|u\|_{H_{V^2}(\Om)}\|v\|_{L^2(\Om)}.$$ By the Hahn-Banach theorem, we can extend $T_0$ to a continuous linear functional on $L^2(\Om)$, and so by Riesz representation theorem there is a unique $M_V^* v\in L^2(\Om)$ such that
$$\la u,V v\ra_{L^2(\Om)}=\la M_V(u),v\ra_{L^2(\Om)}=\la u, M_V^*(v)\ra_{L^2(\Om)}\,\,\text{ for all $u\in H_{V^2}(\Om)$.}$$
Thus we have that $M_V^*(v)=V v\in L^2(\Om)$. This implies that $$\fD(M_V)=H_{V^2}(\Om)=\fD(M_V^*).$$ Thus $M_V$ is self-adjoint and it is also positive, because we have that $$\la M_V(u),u\ra_{L^2(\Om)}=\|u\|^2_{L^2_V(\Om)}\ge 0.$$

(b) The existence and uniqueness of a positive self-adjoint square root operator $P$ of $M_V$ can be obtained by Theorem 13.31 \cite{R}. Hence we are done.  \qed

\,\,\,In what follows, by Lemma 5.1 and Lemma 5.3, for simplicity we may write $$\rL_V=\overline{L_K}+M_V\,\text{ on $\Om$}$$
as an operator $\rL_V:L^2(\Om)\to L^2(\Om)$ with $\fD(\rL_V)=H^{2s}(\BR^n)\cap L^2_{V^2}(\Om)$, which is positive and self-adjoint. Thus there is a positive self-adjoint square root operator $\rL_V^{1/2}$ of $\rL_V$.

\begin{lemma} We have the estimate $$\la\rL_V^{1/2} u,\rL^{1/2}_V v\ra_{L^2(\BR^n)}=\la u,v\ra_{\rX^s_0(\Om)}+\la V u,v\ra_{L^2(\Om)}$$ 
for all $u,v\in C_c^{\iy}(\Om)$. Moreover, if $F:C^{\iy}_c(\Om)\to L^2(\BR^n)$ is the map defined by $F(\vp)=\rL_V^{1/2}\vp$ is injective and satisfies that
\begin{equation}\|F(\phi)\|_{L^2(\BR^n)}=\|\phi\|_{\rY^s_0(\Om)}\,\,\text{ for any $\phi\in C^{\iy}_c(\Om)$.}
\end{equation} 
\end{lemma}

\pf It easily follows from Corollary 5.2, Lemma 5.3 and Theorem 4.1. \qed

\,\,\,Next, we obtain the existence and uniqueness of weak solution of the nonlocal equation $\rL_V=h$ in $\Om$ for the forcing term $h\in L^2(\Om)$ and moreover for $h\in\rY^s_0(\Om)^*$.

\begin{lemma} For each $h\in\rY^s_0(\Om)^*$, there is a unique weak solution $u\in\rY^s_0(\Om)$ of the nonlocal equation $\,\rL_V u=h$ in $\Om$ and $\|u\|_{\rY_0^s(\Om)}\le\|h\|_{\rY_0^s(\Om)^*}$. If $h\in L^2(\Om)$, then we have that
$\|u\|_{\rY^s_0(\Om)}\le C\,\|h\|_{L^2(\Om)}$.
\end{lemma}

\pf We define a bilinear form $\fa:\rY^s_0(\Om)\times\rY^s_0(\Om)\to\BR$ by
\begin{equation*}\fa(u,\phi)=\la\rL_V u,\phi\ra_{L^2(\Om)}. 
\end{equation*}
By Corollary 5.2, it is easy to check that 
\begin{equation*} \fa(u,u)=\|u\|^2_{\rY^s_0(\Om)}\,\,\text{ and }\,\,
|\fa(u,\phi)|\le\|u\|_{\rY^s_0(\Om)}\,\|\phi\|_{\rY^s_0(\Om)}
\end{equation*} for any $u,\phi\in\rY^s_0(\Om)$. Thus the existence result can be obtained by the Lax-Milgram theorem. 

Since $u\in\rY^s_0(\Om)$ is a weak solution of the equation $\,\rL_V u=h$ in $\Om$, we have that
\begin{equation*}\|u\|^2_{\rY_0^s(\Om)}=\la\rL_V u,u\ra_{L^2(\Om)}=\la h,u\ra\le\|h\|_{\rY_0^s(\Om)^*}\|u\|_{\rY_0^s(\Om)},
\end{equation*} and thus we have that $\|u\|_{\rY_0^s(\Om)}\le\|h\|_{\rY_0^s(\Om)^*}$.
If $h\in L^2(\Om)$, then it follows from the dual form of Theorem 4.1 that
$$\|u\|_{\rY_0^s(\Om)}\le\|h\|_{\rY_0^s(\Om)^*}\le C\|h\|_{L^2(\Om)^*}=C\|h\|_{L^2(\Om)}.$$ Hence we are done. \qed

\,\,\,Let $\cD(D)=C^{\iy}_c(D)$ for an open subset $D$ of $\BR^n$, i.e. the class of all smooth functions $\vp$ with compact support in $D$ which inherits a topology induced by the convergence 
\begin{equation*}\vp_k\rightharpoondown\vp\,\,\text{ (as $k\to\iy$)\,\, for $\vp_k,\vp\in\cD(D)$}
\end{equation*}
if and only if there is a compact set $\cQ\subset D$ such that $\supp(\vp_k)\subset\cQ$ for all $k\in\BN$ and $\pa^{\ap}\vp_k$ converges uniformly on $\cQ$ to $\pa^{\ap}\vp$  
for all multi-indices $\ap=(\ap_1,\cdots,\ap_n)$, where $\pa^{\ap}$ denotes
$$\pa^{\ap}=\pa_{x_1}^{\ap_1}\pa_{x_2}^{\ap_2}\cdots\pa_{x_n}^{\ap_n}.$$
A {\it distribution} is defined as a linear functional $T:\cD(D)\to\BR$ which is continuous in the sense that
$$\lim_{k\to\iy}T(\vp_k)=T(\vp)$$ for any sequence $\{\vp_k\}\subset\cD(D)$ with $\vp_k\rightharpoondown\vp\in\cD(D)$.
We denote by $\cD'(D)$ the space of all such distributions on $D$.

Let $T:\cD(D)\to\BR$ be a linear map. Then an equivalent condition of distribution is known as follows (see \cite{R}); $T\in\cD'(D)$ if and only if for each compact set $\cQ\subset D$, there is an integer $N=N(\cQ)>0$ and a constant $C=C(\cQ)>0$ such that
\begin{equation}|T(\vp)|\le C\,\|\vp\|_N\,\,\text{ for all $\vp\in\cD_\cQ$,}
\end{equation}
where $\|\vp\|_N=\max\{\sup_D|\pa^{\ap}\vp|:|\ap|:=\sum_{i=1}^n\ap_i\le N\}$ and $\cD_\cQ$ denotes the class of all smooth functions in $D$ which is supported in $\cQ$. In what follows, we write $\la T,\vp\ra=T(\vp)$.

\begin{lemma} If $u\in\dot H^s(\BR^n)$ for $s\in(0,1)$, then the linear map $\rL_V u:\cD(\BR^n)\to\BR$ defined by 
\begin{equation}\la\rL_V u,\vp\ra=\int_{\BR^n}u(y)\,\rL_V\vp(y)\,dy
\end{equation} is a distribution, i.e. $\rL_V u\in\cD'(\BR^n)$.
\end{lemma}

\pf Since $0<s<\f{n}{2}$ for $n\ge 2$, it is easy to check that 
\begin{equation}\rL_V\vp\in L^{\f{2n}{n+2s}}(\BR^n)
\end{equation} for any $\vp\in \cD(\BR^n)$. Thus, by Proposition 2.2 and H\"older's inequality, the linear functional $\rL_V u$ is well-defined on $\cD(\BR^n)$. 

Now it remains only to show that $\rL_V u$ satisfies the continuity property (5.7). Indeed, taking arbitrary compact set $\cQ\subset\BR^n$ and $\vp\in\cD_\cQ$, it follows from (2.7) and the mean value theorem that
\begin{equation*}\begin{split}
&|\la\rL_V u,\vp\ra|^2\\&\le c^2_{n,s}\Ld^2\biggl(\iint_{\BR^n\times\BR^n}\f{(u(x)-u(y))(\vp(x)-\vp(y))}{|x-y|^{n+2s}}\,dx\,dy\biggr)^2\\
&\le C^2_0\iint_{\BR^{2n}_\cQ}\f{|\vp(x)-\vp(y)|^2}{|x-y|^{n+2s}}\,dx\,dy\\
&\le C^2_0\biggl(\iint_{B_{2 r_\cQ}\times B_{2 r_\cQ}}\f{|\vp(x)-\vp(y)|^2}{|x-y|^{n+2s}}\,dx\,dy+\iint_{B^c_{2 r_\cQ}\times B_{r_\cQ}}\f{|\vp(y)|^2}{|x-y|^{n+2s}}\,dx\,dy\biggr)\\
&\le C_0^2\biggl(\int_{B_{2 r_\cQ}}\int_{B_{4 r_\cQ}}\f{\sup_{\BR^n}|\n\vp|^2}{|x|^{n+2(s-1)}}\,dx\,dy+\int_{B_{r_\cQ}}\int_{B^c_{2 r_\cQ}}\f{\sup_{\BR^n}|\vp|^2}{|x|^{n+2s}}\,dx\,dy\biggr)\\
&\le C_1 C^2_0\,\|\vp\|^2_1,
\end{split}\end{equation*}
where $r_\cQ=\sup_{y\in\cQ}|y|$, $C_1=C_1(n,s,\cQ)>0$ is a constant and 
$$C_0=\sqrt{\f{c^3_{n,s}}{2}}\,\Ld\,\|u\|_{\dot H^s(\BR^n)}.$$
Therefore we complete the proof. \qed

\begin{lemma} Let $V\in L^q_{\loc}(\BR^n)$ be a nonnegative potential with $q>\f{n}{2s},$ $s\in(0,1)$. If $\,h\in L^{\f{2n}{n+2s}}(\BR^n)$, then there exists a function $u\in\dot H^s(\BR^n)$ such that 
\begin{equation*}\rL_V u=h\,\,\,\text{ in the sense of $\cD'(\BR^n)$.}
\end{equation*} Moreover, there is an increasing sequence $\{a_k\}_{k\in\BN}$ with $\lim_{k\to\iy}a_k=\iy$ such that $u_{a_k}\in\rY^s_0(B_{a_k})$ satisfies the nonlocal equation $\rL_V u_{a_k}=h$ in $B_{a_k}$ in the weak sense and $\,\lim_{k\to\iy}u_{a_k}=u\,\,\,\aee$ in $\BR^n$.
\end{lemma}

\pf We observe that 
\begin{equation}h\in L^{\f{2n}{n+2s}}(B_k)\subset L^2(B_k)\subset\rY^s_0(B_k)^*\,\,\,\text{ for all $k\in\BN$.}
\end{equation}
By Theorem 4.2 and Lemma 5.5, for each $k\in\BN$ there exists a nonnegative function $u_k\in\rY^s_0(B_k)$ which is a weak solution of the equation $\rL_V u_k=h$ in $B_k$, i.e. 
\begin{equation}\la u_k,\vp\ra_{\rX^s_0(B_k)}+\la V u_k,\vp\ra_{L^2(B_k)}=\la h,\vp\ra_{L^2(B_k)}
\end{equation} for any $\vp\in\rY^s_0(B_k)$.
Taking $\vp=u_k$ in (5.11), it follows from H\"older's inequality, Lemma 2.1 and the fractional Sobolev embedding on $\dot H^s(\BR^n)$ (Proposition 2.2) that
\begin{equation}\begin{split}\f{c_{n,s}}{2}\,\|u_k\|^2_{\dot H^s(\BR^n)}&=\|u_k\|^2_{\rX^s_0(B_k)}\le\|u_k\|^2_{\rY^s_0(B_k)}\\
&\le\|h\|_{L^{\f{2n}{n+2s}}(\BR^n)}\|u_k\|_{L^{\f{2n}{n-2s}}(\BR^n)}\\
&\le C\,\|h\|_{L^{\f{2n}{n+2s}}(\BR^n)}\|u_k\|_{\dot H^s(\BR^n)}
\end{split}\end{equation} where $C>0$ is a constant depending only on $n,s$, but not on $k$. This implies that
\begin{equation}\sup_{k\in\BN}\|u_k\|_{L^{\f{2n}{n-2s}}(\BR^n)}\le C\,\sup_{k\in\BN}\|u_k\|^2_{\dot H^s(\BR^n)}\le C\,\|h\|_{L^{\f{2n}{n+2s}}(\BR^n)}.
\end{equation}
By weak compactness, there are a subsequence $\{u_{k_i}\}_{i\in\BN}$ and $u\in\dot H^s(\BR^n)$ such that
\begin{equation}u_{k_i}\rightharpoonup u\,\,\text{ in $L^{\f{2n}{n-2s}}(\BR^n)$ }
\,\,\text{ and }\,\,u_{k_i}\rightharpoonup u\,\,\text{ in $\dot H^s(\BR^n)$ }
\end{equation}
Let us take any $\vp\in\cD(\BR^n)$. Then there is some $m\in\BN$ such that $\vp\in\cD(B_m)$. We note that $\rL_V u_k=h$ in the sense of $\cD'(B_m)$ for any $k\ge m$. Thus, by (5.14), Corollary 5.2 and Lemma 5.6, we conclude that
\begin{equation}\begin{split}
\la\rL_V u,\vp\ra&=\int_{\BR^n}u(y)\,\rL_V\vp(y)\,dy=\lim_{k\to\iy}\int_{\BR^n}u_k(y)\,\rL_V\vp(y)\,dy=\la h,\vp\ra.
\end{split}\end{equation}
Thus this implies that $\rL_V u=h$ in the sense of $\cD'(\BR^n)$. 

Finally, we see from (5.14) and Theorem 7.1 \cite{DPV} that the sequence $\{u_k\}$ is precompact in $L^2(B)$ for any ball $B$ in $\BR^n$. This implies the required almost everywhere convergence. \qed

\,\,\,Next, we shall prove our main theorem which is {\it an extension of the Malgrange-Ehrenpreis theorem} to nonlocal Schr\"odinger operators $\rL_V$.

For $l\in\BN$, let $f_l(x)=l^n f\bigl(l x\bigr)$ where $f\in C_c(B_1)$ is a nonnegative function with $\|f\|_{L^1(\BR^n)}=1$.
From Lemma 5.7, we see that there is a sequence $\{u_l\}_{l\in\BN}\subset\dot H^s(\BR^n)$ such that
\begin{equation}\rL_Vu_l=f_l\,\,\,\,\text{ in the sense of $\cD'(\BR^n)$. } 
\end{equation}
As a matter of fact, we see from a weak maximum principle (Theorem 4.2) and the proof of Lemma 5.7 that each $u_l$ is nonnegative and, for each $l\in\BN$, there are a sequence $\{a^l_k\}_{k\in\BN}$ with $\lim_{k\to\iy}a^l_k=\iy$ and a sequence $\{u_l^k\}_{k\in\BN}$ of nonnegative functions $u_l^k\in\rY^s_0(B_{a^l_k})$ such that 
\begin{equation}\begin{split}&\quad u_l=\ds\lim_{k\to\iy}u_l^k\,\,\text{ $\aee$ in $\BR^n$,}\\ 
&\quad\rL_V u_l^k=f_l\,\text{ in $B_{a^l_k}\,$ in the weak sense.}
\end{split}\end{equation}
Then we obtain uniform estimates for the sequence $\{u_l\}_{l\in\BN}$ in the following lemma in order to prove the main theorem. 

\begin{lemma} $(a)$ For each $l\in\BN$ and $p\in[1,\f{n}{n-2s})$ with $s\in(0,1)$, there exists some constant $C=C(p,\ld,n,s)>0$ such that 
\begin{equation}\|u_l\|_{L^p(B_r(x_0))}\le C\,r^{\f{n}{p}-(n-2s)}\,\,\,\text{for any $x_0\in\BR^n$ and $r>0$.}
\end{equation} 

$(b)$ For each $l\in\BN$ and $p\in[1,\f{n}{n-s})$ and $\gm\in(0,s)$ with $s\in(0,1)$, there exists some constant $C=C(p,n,s)>0$ such that
\begin{equation}[u_l]_{W^{\gm,p}(B_r(x_0))}\le C\,r^{\f{n}{p}-(n-2s+\gm)}\,\,\,\text{for any $x_0\in\BR^n$ and $r>0$.}
\end{equation} 

$(c)$ If $\,V\in L^q_{\loc}(\BR^n)$ is nonnegative for $q>\f{n}{2s}$ with $s\in(0,1)$, then for each $l\in\BN$ there is a constant $C=C(p,\ld,n,s)>0$ such that
\begin{equation}\|u_l\|_{L^1_V(B_r(x_0))}\le\f{C\,\|V\|_{L^q(B_r(x_0))}}{r^{n-2s-n/p}}\,\,\,\text{for any $x_0\in\BR^n$ and $r>0$,} 
\end{equation} where $p\in[1,\f{n}{n-2s})$ is the dual exponent of $q>\f{n}{2s}$.

$(d)$ For each $l\in\BN$ and for any $s\in(0,1)$ and $x\in\BR^n\s\{0\}$ with $|x|\ge 3/l$, we have that
\begin{equation}\bigl|u_l(x)\bigr|\le\f{C}{|x|^{n-2s}}
\end{equation} where $C=C(n,s,\ld,\Ld)>0$ is a constant depending only on $n,s,\ld$ and $\Ld$.
\end{lemma}

\pf (a) Motivated by \cite{KMS}, we consider a nonnegative function $\vp:\BR^n\to\BR$ defined by 
\begin{equation*}\vp=\bt^{1-\ap}-(\bt+u^k_l)^{1-\ap}
\end{equation*}
where $\bt>0$ and $\ap\in(1,2)$ will be chosen later. Since each $u^k_l$ is supported in $B_{a_k^l}$ and $a(t)=(\bt+t)^{1-\ap}$ is Lipschitz continuous on $(0,\iy)$, by Lemma 2.1 and the definition of $\rY^s_0(B_{a_k^l})$ we see that $\vp\in \rY^s_0(B_{a_k^l})\subset H^s(\BR^n)$. Thus we can use $\vp$ as a testing function. From the weak formulation of the nonlocal equation given in (5.17), we have that
\begin{equation}\begin{split}&\iint_{\BR^n\times\BR^n}(u^k_l(x)-u^k_l(y))(\vp(x)-\vp(y))K(x-y)\,dx\,dy\le\bt^{1-\ap}.
\end{split}\end{equation}
By the mean value theorem, we have that
\begin{equation}\begin{split}\vp(x)-\vp(y)=(u^k_l(x)-u^k_l(y))\int_0^1\f{(\ap-1)\,dt}{(\bt+t u^k_l(y)+(1-t)u^k_l(x))^\ap}.
\end{split}\end{equation}
Thus it follows from (5.22) and (5.23) that
\begin{equation}\bt^{1-\ap}\ge (\ap-1)\iint_{\BR^n\times\BR^n}\biggl(\int_0^1\f{(u^k_l(x)-u^k_l(y))^2\,dt}{(\bt+t u^k_l(x)+(1-t) u^k_l(y))^{\ap}}\biggr)\f{dx\,dy}{|x-y|^{n+2s}}.
\end{equation}
By Jensen's inequality and the mean value theorem, we have that
\begin{equation}\begin{split}
&\biggl(\f{u^k_l(x)}{(\bt+u^k_l(x))^{\ap/2}}-\f{u^k_l(y)}{(\bt+u^k_l(y))^{\ap/2}}\biggr)^2\\
&\qquad\qquad\qquad=\biggl(\int_0^1\f{d}{dt}\biggl[\f{t u^k_l(x)+(1-t)u^k_l(y)}{(\bt+t u^k_l(x)+(1-t)u^k_l(y))^{\ap/2}}\biggr]\,dt\biggr)^2\\
&\qquad\qquad\qquad\le 4\int_0^1\f{|u^k_l(x)-u^k_l(y)|^2\,dt}{(\bt+t u^k_l(x)+(1-t)u^k_l(y))^{\ap}}.
\end{split}\end{equation}
Plugging (5.25) into (5.24) and using Proposition 2.2 (fractional Sobolev embedding theorem), we obtain that
\begin{equation}\begin{split}
\bt^{1-\ap}&\ge\f{\ap-1}{4}\iint_{\BR^n\times\BR^n}\biggl(\f{u^k_l(x)}{(\bt+u^k_l(x))^{\ap/2}}-\f{u^k_l(y)}{(\bt+u^k_l(y))^{\ap/2}}\biggr)^2\f{dx\,dy}{|x-y|^{n+2s}}\\
&\ge\f{\ap-1}{4}\biggl(\int_{\BR^n}\biggl(\f{u^k_l(x)}{(\bt+u^k_l(x))^{\ap/2}}\biggr)^{\f{2n}{n-2s}}\,dx\biggr)^{\f{n-2s}{n}}.
\end{split}\end{equation}

Take any $p\in[1,\f{n}{n-2s})$ and choose some $\ap\in(1,2)$ satisfying $\f{\ap}{2}=1-\f{(n-2s)p}{2n}$. By H\"older's inequality and (5.26), we have that
\begin{equation}\begin{split}
&\|u^k_l\|^p_{L^p(B_r(x_0))}\\
&\le\biggl(\int_{\BR^n}\biggl(\f{u^k_l(x)}{(\bt+u^k_l(x))^{\ap/2}}\biggr)^{\f{2n}{n-2s}}\biggr)^{\f{(n-2s)p}{2n}}\biggl(\int_{B_r(x_0)}(\bt+u^k_l(x))^{\f{\ap}{2}p q'}\,dx\biggr)^{\f{1}{q'}}\\
&\le C\,\bt^{(1-\ap)\f{p}{2}}\biggl(\int_{B_r(x_0)}(\bt+u^k_l(x))^{\f{\ap}{2}p q'}\,dx\biggr)^{\f{1}{q'}}
\end{split}\end{equation}
where $q'$ is the dual exponet of $q=\f{2n}{(n-2s)p}$, because $\f{\ap}{2}q'=1$.
If we set
\begin{equation*}\bt=\biggl(\f{1}{|B_r(x_0)|}\int_{B_r(x_0)}u^k_l(x)^p\,dx\biggr)^{\f{1}{p}}=\f{|B_1(0)|^{-1/p}}{r^{n/p}}\biggl(\int_{B_r(x_0)}u^k_l(x)^p\,dx\biggr)^{\f{1}{p}},
\end{equation*}
then by (5.29) we have that
$$\|u^k_l\|^p_{L^p(B_r(x_0))}\le\f{C}{r^{\f{n}{2}(1-\ap)}}\|u^k_l\|^{p/2}_{L^p(B_r(x_0))}.$$ This gives that, for each $l\in\BN$,
\begin{equation}\|u^k_l\|_{L^p(B_r(x_0))}\le C\,r^{\f{n}{p}-(n-2s)}\,\,\text{ for any $k\in\BN$.}
\end{equation} 
Hence it follows from (5.17), (5.28) and Fatou's lemma that
\begin{equation*}\|u_l\|_{L^p(B_r(x_0))}\le C\,r^{\f{n}{p}-(n-2s)}\,\,\text{ for each $l\in\BN$.}
\end{equation*}

(b) For $\gm\in(0,s]$ with $s\in(0,1)$, we set
\begin{equation*}U_{\gm}(x,y)=\f{|u^k_l(x)-u^k_l(y)|}{|x-y|^{\gm}}.
\end{equation*}
For $p\in[1,2)$, we write 
\begin{equation*}
U^p_s(x,y)=\biggl(\f{U^2_s(x,y)}{(\bt+u^k_l(x)+u^k_l(y))^{\ap}}\biggr)^{\f{p}{2}}\,\bigl((\bt+u^k_l(x)+u^k_l(y))^{\ap}\,|x-y|^{2(s-\gm)}\bigr)^{\f{p}{2}}.
\end{equation*}
Then it follows from H\"older's inequality and (5.24) that
\begin{equation}\begin{split}
&\iint_{B^2_r(x_0)}U^p_{\gm}(x,y)\,\f{dx\,dy}{|x-y|^n}\le\biggl(\iint_{B^2_r(x_0)}\f{U^2_s(x,y)}{(\bt+u^k_l(x)+u^k_l(y))^{\ap}}\,\f{dx\,dy}{|x-y|^n}\biggr)^{\f{p}{2}}\\
&\qquad\qquad\qquad\qquad\quad\times\biggl(\iint_{B^2_r(x_0)}\f{(\bt+u^k_l(x)+u^k_l(y))^{\ap p/(2-p)}}{|x-y|^{n-2(s-\gm)p/(2-p)}}\,dx\,dy\biggr)^{\f{2-p}{2}}\\
&\qquad\qquad\qquad\qquad\quad\le C\,\bt^{(1-\ap)\f{p}{2}}r^{(s-\gm)p}\biggl(\int_{B_r(x_0)}(\bt+u^k_l(x))^{\f{\ap p}{2-p}}\,dx\biggr)^{\f{2-p}{2}},
\end{split}\end{equation} where $B^2_r(x_0)=B_r(x_0)\times B_r(x_0)$. If we choose $\bt$  satisfying
\begin{equation*}\bt=\biggl(\f{1}{|B_r(x_0)|}\int_{B_r(x_0)}u^k_l(x)^{\f{\ap p}{2-p}}\,dx\biggr)^{\f{2-p}{\ap p}},
\end{equation*}
then (5.29) leads us to get that
\begin{equation}\begin{split}
&\iint_{B^2_r(x_0)}U^p_{\gm}(x,y)\,\f{dx\,dy}{|x-y|^n}\\
&\le C\,r^{(s-\gm)p}\biggl(\f{1}{r^n}\int_{B_r(x_0)}u^k_l(x)^{\f{\ap p}{2-p}}\,dx\biggr)^{\f{(2-p)(1-\ap)}{2\ap}}\biggl(\int_{B_r(x_0)}u^k_l(x)^{\f{\ap p}{2-p}}\,dx\biggr)^{\f{2-p}{2}}.
\end{split}\end{equation}
Take any $p\in[1,\f{n}{n-s})$. We may select some $\ap\in(1,2)$ satisfying $\f{\ap p}{2-p}\in[1,\f{n}{n-2s})$, because $h(t)=\f{t}{2-t}\ge t$ for any $t\in[1,2)$ and $h(\f{n}{n-s})=\f{n}{n-2s}$. Applying (5.28) to (5.30), we obtain that
\begin{equation*}\iint_{B^2_r(x_0)}U^p_{\gm}(x,y)\,\f{dx\,dy}{|x-y|^n}\le C\,r^M
\end{equation*}
where the power index $M$ of $r$ is given by
\begin{equation*}\begin{split}M&=(s-\gm)p-\f{n(2-p)(1-\ap)}{2\ap}\\
&+\biggl[\f{n(2-p)}{\ap p}-(n-2s)\biggr]\f{\ap p}{2-p}\biggl[\f{(2-p)(1-\ap)}{2\ap}+\f{2-p}{2}\biggr]\\
&=n-(n-2s)p-\gm p.
\end{split}\end{equation*}
Thus we have that, for each $l\in\BN$,
\begin{equation}[u^k_l]_{W^{\gm,p}(B_r(x_0))}\le C\,r^M\,\,\text{ for any $k\in\BN$.}
\end{equation}
Hence, (5.17), (5.31) and Fatou's lemma imply the second part.

(c) By H\"older's inequality, we have that
\begin{equation}\begin{split}
\|u_l\|_{L^1_V(B_r(x_0))}&\le\|u_l\|_{L^p(B_r(x_0)}\|V\|_{L^q(B_r(x_0))}
\end{split}\end{equation}
where $p\in[1,\f{n}{n-2s})$ is the dual exponent of $q>\f{n}{2s}$. Therefore, by (5.18) and (5.32), we obtain that
$$\|u_l\|_{L^1_V(B_r(x_0))}\le\f{C\,\|V\|_{L^q(B_r(x_0))}}{r^{n-2s-n/p}}$$ with a universal constant $C=C(p,\ld,n,s)>0$. Thus we are done. 

(d) Take any $x\in\BR^n$ with $|x|\ge 3/l$. Then we note that $B_{2r}(x)\cap B_{1/l}=\phi$ where $r=\f{1}{4}|x|$. By Lemma 4.5 \cite{CK}, we have the estimate
\begin{equation}\|u^k_l\|_{L^{\iy}(B_r(x))}\le\f{C}{r^{n/2}}\,\|u^k_l\|_{L^2(B_{2r}(x))}
\end{equation} with some universal constant $C>0$. Using a standard argument in \cite{HL}, for any $p\in[1,\iy)$ there is a universal constant $C=C(n,s,p)>0$ such that
\begin{equation}\|u^k_l\|_{L^{\iy}(B_r(x))}\le\f{C}{r^{n/p}}\,\|u^k_l\|_{L^p(B_{2r}(x))}.
\end{equation}
Thus it follows from (5.28) and (5.34) with $p=1$ that
\begin{equation*}\|u^k_l\|_{L^{\iy}(B_r(x))}\le\f{C}{r^n}\,\|u^k_l\|_{L^1(B_{2r}(x))}\le\f{C}{r^{n-2s}}.
\end{equation*} Therefore this and (5.17) imply (5.21). Hence we complete the proof.
\qed

\,\,\,{\bf $[$Proof of Theorem 1.1\,$]$\,\,\,} In what follows, we always take some $\gm\in(0,s)$ (sufficiently close to $s$) with 
\begin{equation}1-s>s-\gm
\end{equation} and a sufficiently small $\e\in(0,1)$ with
\begin{equation}n\e+(2s-\gm-1)(1+\e)\le 0\,\,\text{ and }\,\,p:=1+\e<\f{n}{n-s}.
\end{equation}
From (5.17), Lemma 5.8 [(5.18), (5.19),(5.20)] and Theorem 7.1 \cite{DPV}, the standard diagonalization process yields that there exist a subsequence $\{u_{l_i}\}_{i\in\BN}$ of $\{u_l\}_{l\in\BN}$ and $\fe_V\in W^{\gm,p}_{\loc}(\BR^n)\subset\cD'(\BR^n)$ such that 
\begin{equation}\begin{cases}
u_{l_i}\rightharpoonup\fe_V &\text{ in $W^{\gm,p}(B_r(x_0))$, }\\
u_{l_i}\rightharpoonup\fe_V &\text{ in $L^p(B_r(x_0))$, }\\
V u_{l_i}\rightharpoonup V\fe_V &\text{ in $L^1(B_r(x_0))$, }\\
u_{l_i}\to\fe_V\,\,\,\,\aee \,\, &\text{ in $\,B_r(x_0)$}
\end{cases}\end{equation} for any $(x_0,r)\in\BR^n\times(0,\iy)$.
We write $\{u_l\}_{l\in\BN}$ instead of $\{u_{l_i}\}_{i\in\BN}$, for simplicity.

First, we show that the linear map $\rL_V\fe_V:\cD(\BR^n)\to\BR$ defined by
\begin{equation}\la\rL_V\fe_V,\vp\ra=\int_{\BR^n}\fe_V(y)\,\rL_V\vp(y)\,dy
\end{equation}
is a distribution, i.e. $\rL_V\fe_V\in\cD'(\BR^n)$. Indeed, take any compact set $\cQ\subset\BR^n$ and $\vp\in\cD_\cQ$ for this proof. If $1<p<\f{n}{n-s}$, then it follows from H\"older's inequality and Young's inequality that 
\begin{equation}\begin{split}
|\la\rL_V\fe_V,\vp\ra|&\le c_{n,s}\Ld\iint_{\BR^{2n}_\cQ}\f{|(\fe_V(x)-\fe_V(y))(\vp(x)-\vp(y))|}{|x-y|^{n+2s}}\,dx\,dy\\
&\le c_{n,s}\Ld\biggl(\iint_{\BR^{2n}_\cQ}\f{|\fe_V(x)-\fe_V(y)|^p}{|x-y|^{n+p\gm }}\,dx\,dy\biggr)^{1/p}\\
&\qquad\times\biggl(\iint_{\BR^{2n}_\cQ}\f{|\vp(x)-\vp(y)|^{p'}}{|x-y|^{n+p\gm+p'(2s-p\gm)}}\,dx\,dy\biggr)^{1/p'}\\
&:= c_{n,s}\Ld\,\,I(\cQ,p)\,\,J(\cQ,p')
\end{split}\end{equation} where $p'$ is the dual exponent of $p$,
because $1<p<\f{n}{n-s}<2$ for any $s\in(0,1)$ and $n\ge 2$, and
\begin{equation*}n+2s=\f{n+p\gm}{p}+\f{n+p\gm}{p'}+(2s-p\gm).
\end{equation*}
By simple calculation, we have that
\begin{equation}\begin{split}
I^p(\cQ,p)&\le\iint_{B_{2 r_\cQ}\times B_{2 r_\cQ}}\f{|\fe_V(x)-\fe_V(y)|^p}{|x-y|^{n+p\gm }}\,dx\,dy\\
&\qquad+2^{p+1}\int_{B_{r_\cQ}}\int_{B^c_{2 r_\cQ}}\f{|\fe_V(y)|^p}{|x-y|^{n+p\gm}}\,dx\,dy\\
&\le[\fe_V]^p_{W^{\gm,p}(B_{2 r_\cQ})}+C(n,p,\gm,\cQ)\,\|\fe_V\|^p_{L^p(B_{r_\cQ})}\\
&\le C(n,p,\gm,\cQ)\,\|\fe_V\|^p_{W^{\gm,p}(B_{2 r_\cQ})},
\end{split}\end{equation} where $r_\cQ=\sup_{y\in\cQ}|y|$.
By the mean value theorem, we obtain that
\begin{equation}\begin{split}
J^{p'}(\cQ,p')&\le\iint_{B_{2 r_\cQ}\times B_{2 r_\cQ}}\f{|\vp(x)-\vp(y)|^{p'}}{|x-y|^{n+p\gm+p'(2s-p\gm)}}\,dx\,dy\\
&\quad+\iint_{B^c_{2 r_\cQ}\times B_{r_\cQ}}\f{|\vp(y)|^{p'}}{|x-y|^{n+p\gm+p'(2s-p\gm)}}\,dx\,dy\\
&\le\int_{B_{2 r_\cQ}}\int_{B_{4 r_\cQ}}\f{\sup_{\BR^n}|\n\vp|^{p'}}{|x|^{n+p\gm+p'(2s-p\gm-1)}}\,dx\,dy\\
&\quad+\int_{B_{r_\cQ}}\int_{B^c_{2 r_\cQ}}\f{\sup_{\BR^n}|\vp|^{p'}}{|x|^{n+p\gm+p'(2s-p\gm)}}\,dx\,dy\\
&\le C(n,s,p,\gm,\cQ)\,\|\vp\|_1^{p'},
\end{split}\end{equation}
because by (5.35) we have that
\begin{equation*}
p\gm+p'(2s-p\gm-1)=\f{p(2s-\gm-1)}{p-1}<0\text{ and }p\gm+p'(2s-p\gm)=\f{(2s-\gm)p}{p-1}>0.
\end{equation*} From (5.39), (5.40) and (5.41), we conclude that
$$|\la\rL_V\fe_V,\vp\ra|\le C(n,s,p,\gm,\cQ)\,\|\fe_V\|_{W^{\gm,p}(B_{2 r_\cQ})}\,\|\vp\|_1.$$
Hence this implies that $\rL_V\fe_V\in\cD'(\BR^n)$.

Take any $\vp\in\cD(\BR^n)$. Then there exists an integrer $m\in\BN$ such that $\vp\in\cD(B_m)$.
Using Lemma 5.6, as in (5.15) we have that 
\begin{equation*}\begin{split}
\la\rL_V u_l,\vp\ra&=\int_{\BR^n}u_l(x)\,\rL_V\vp(x)\,dx=\lim_{k\to\iy}\int_{\BR^n}u_l^k(x)\,\rL_V\vp(x)\,dx\\&=\lim_{k\to\iy}\int_{\BR^n}\rL_V u_l^k(x)\,\vp(x)\,dx=\la f_l,\vp\ra
\end{split}\end{equation*} for each $l\in\BN$.
Thus, for each $l\in\BN$, we see that
\begin{equation}\begin{split}&\la\rL_V u_l,\vp\ra=\iint_{\BR^{2n}_{B_m}}(u_l(x)-u_l(y))(\vp(x)-\vp(y))K(x-y)\,dx\,dy\\
&\qquad\qquad\qquad+\int_{B_m}V(y)u_l(y)\vp(y)\,dz=\int_{\BR^n}f_l(y)\vp(y)\,dy.
\end{split}\end{equation} 
Taking the limit $l\to\iy$ on (5.42), by (5.38) we then claim that
\begin{equation*}\begin{split}&\la\rL_V\fe_V,\vp\ra=\iint_{\BR^{2n}_{B_m}}(\fe_V(x)-\fe_V(y))(\vp(x)-\vp(y))K(x-y)\,dx\,dy\\
&\qquad\qquad\qquad\qquad\qquad+\int_{B_m}V(y)\fe_V(y)\vp(y)\,dy=\vp(0)=\la\dt_0,\vp\ra.
\end{split}\end{equation*}
For our aim, we write 
\begin{equation}\begin{split}
&\la\rL_V u_l,\vp\ra=\iint_{\BR^{2n}_{B_m}}(\fe_V(x)-\fe_V(y))(\vp(x)-\vp(y))K(x-y)\,dx\,dy\\
&\qquad\qquad\qquad\qquad\qquad+\int_{B_m}V(y)\fe_V(y)\vp(y)\,dy\\
&\quad+\iint_{\BR^{2n}_{B_m}}(u_l(x)-\fe_V(x)-u_l(y)+\fe_V(y))(\vp(x)-\vp(y))K(x-y)\,dx\,dy\\
&\qquad\qquad\qquad\qquad\qquad+\int_{B_m}V(y)(u_l(y)-\fe_V(y))\vp(y)\,dy.
\end{split}\end{equation}
Here we denote by $\cI_l$ the last third integral and $\cJ_l$ the last fourth integral in the right side of (5.43). From Lemma 5.8 and Fatou's lemma,  we note that
\begin{equation}\begin{split}[\fe_V]_{W^{\gm,p}(B_k)}+[u_l]_{W^{\gm,p}(B_k)}\le C_1,\\
\|\fe_V\|_{L^p(B_k)}+\|u_l\|_{L^p(B_k)}\le C_2,\\
\|\fe_V\|_{L^1_V(B_k)}+\|u_l\|_{L^1_V(B_k)}\le C_3\,\|V\|_{L^1(B_k)},
\end{split}\end{equation}
where $C_1,C_2,C_3>0$ are some constants depending only on $p,\ld,\Ld,n,s,r$ as in (5.18), (5.19) and (5.20). From the weak convergence $V u_{k_j}\rightharpoonup V\fe_V$  in $L^1(B_m)$ obtained in (5.37), we easily derive that $\lim_{l\to\iy}\cJ_l=0$. So it remains only to show that $\lim_{l\to\iy}\cI_l=0$. To execute this, we split $\cI_l$ into
\begin{equation}\begin{split}
\cI_l&=\iint_{B^2_m}(u_l(x)-\fe_V(x)-u_l(y)+\fe_V(y))(\vp(x)-\vp(y))K(x-y)\,dx\,dy\\
&\qquad\qquad\qquad+\int_{B_R\s B_m}\int_{B_m}(u_l(x)-\fe_V(x))\,\vp(x)\,K(x-y)\,dx\,dy\\
&\qquad\qquad\qquad+\int_{\BR^n\s B_R}\int_{B_m}(u_l(x)-\fe_V(x))\,\vp(x)\,K(x-y)\,dx\,dy\\
&\qquad\qquad\qquad+\int_{B_m}\int_{B_R\s B_m}(u_l(y)-\fe_V(y))\,\vp(y)\,K(x-y)\,dx\,dy\\
&\qquad\qquad\qquad+\int_{B_m}\int_{\BR^n\s B_R}(u_l(y)-\fe_V(y))\,\vp(y)\,K(x-y)\,dx\,dy\\
&:=\cI^1_l+\cI^2_l+\cI^3_l+\cI^4_l+\cI^5_l=\cI^1_l+2 \cI^2_l+2 \cI^3_l
\end{split}\end{equation}
where $B_m\subset B_R$ for a sufficiently large $R>0$ to be determined later.

We have only to show that $\cI^1_l,\cI^2_l,\cI^3_l\to 0$ as $l\to\iy$. We now prove that  $\cI^1_l\to 0$ as $l\to\iy$. To do this, we shall apply Vitali convergence theorem. So we need only to show that the functions $v_l(x,y)$ given by $$v_l(x,y)=(u_l(x)-\fe_V(x)-u_l(y)+\fe_V(y))(\vp(x)-\vp(y))K(x-y)$$
are equibounded in $L^1(B^2_m;dx\,dy)$ and equi-integrable in $B^2_m$. For this, it is enough to show that the sequence $\{v_l\}$ 
is equibounded in $L^{1+\e}(B^2_m;dx\,dy)$ for some $\e>0$. By (5.35) and (5.36), we have that
\begin{equation}\begin{split}
[v_l(x,y)]^{1+\e}&\le c_{n,s}\,\Ld\,\f{[(|\fe_V(x)-\fe_V(y)|+|u_l(x)-u_l(y)|)(\vp(x)-\vp(y))]^{1+\e}}{|x-y|^{(n+2s)(1+\e)}}\\
&\le c_{n,s}\,\Ld\,\|\n\vp\|^{1+\e}_{L^{\iy}(\BR^n)}\,\f{(|\fe_V(x)-\fe_V(y)|+|u_l(x)-u_l(y)|)^{1+\e}}{|x-y|^{n+\gm(1+\e)}|x-y|^{n\e+(2s-\gm-1)(1+\e)}}\\
&\le C\,\biggl(\f{|\fe_V(x)-\fe_V(y)|^{1+\e}}{|x-y|^{n+\gm(1+\e)}}+\f{|u_l(x)-u_l(y)|^{1+\e}}{|x-y|^{n+\gm(1+\e)}}\biggr)
\end{split}\end{equation}
where $C=c(n,m,\gm,\e,s)\,2^{1+\e}\,c_{n,s}\,\Ld\|\n\vp\|^{1+\e}_{L^{\iy}(\BR^n)}>0$. From (5.44) and (5.46), we thus conclude that $\{v_l\}_{l\in\BN}$ is equibounded in $L^{1+\e}(B^2_m;dx\,dy)$. 

For the convergence $\cI^2_l\to 0$ as $l\to\iy$, we need only to show again that the functions $w_l(x,y)$ given by $$w_l(x,y)=(u_l(x)-\fe_V(x))\,\vp(x)\,K(x-y)$$ are equibounded in $L^{1+\e}((B_R\s B_m)\times B_m;dx\,dy)$. Indeed, if we denote by $d=\dist(\supp(\vp),\pa B_m)>0$ (since $\supp(\vp)\subset B_m$ is compact), then by (5.44) we have that
\begin{equation*}\begin{split}&\iint_{(B_R\s B_m)\times B_m}[w_l(x,y)]^{1+\e}\,dx\,dy\\
&\qquad\qquad\qquad\qquad\le\f{2^{1+\e}\,c_{n,s}\,\Ld}{d^{n+2s}}\,\|\fe_V+u_l\|^{1+\e}_{L^{1+\e}(B_m)}\,|B_R\s B_m|\\
&\qquad\qquad\qquad\qquad\le\f{(2 C_2)^{1+\e}\,c_{n,s}\,\Ld}{d^{n+2s}}\,|B_R\s B_m|<\iy\,\,\text{ for any $R>0$.}
\end{split}\end{equation*}

Finally, we have only to prove that $\cI^3_l\to 0$ as $l\to\iy$; in fact, we shall prove that for any $\vep>0$ there exists some $R=R(\vep)>0$ so large that $\limsup_{l\to\iy}|\cI^3_l|<\vep$. Fix any $\vep>0$ and take $R>0$ so that $R>2m$. Then we see that $|y|>2|x|$ for any $y\in\BR^n\s B_R$ and $x\in B_m$.
Applying H\"older's inequality, by (5.44) we obtain that
\begin{equation*}\begin{split}
|\cI^3_l|&\le c_{n,s}\,\Ld\iint_{(\BR^n\s B_R)\times B_m}\f{|\fe_V(x)|+|u_l(x)|}{|x-y|^{n+2s}}\,dx\,dy\\
&\le\biggl(\int_{\BR^n\s B_R}\f{2^{n+2s}\,c_{n,s}\,\Ld}{|y|^{n+2s}}\,dy\biggr)
|B_m|^{\f{\e}{1+\e}}\bigl(\,\|\fe_V\|_{L^{1+\e}(B_m)}+\|u_l\|_{L^{1+\e}(B_m)}\,\bigr)\\
&\le\f{2^{n+2s}C_2|S^{n-1}|\,c_{n,s}\,\Ld\,|B_m|^{\f{\e}{1+\e}}}{2s\,R^{2s}}<\vep,
\end{split}\end{equation*}
whenever $R>2m$ is taken sufficiently large. 

Also, the second part (1.5) can be easily derived from (d) of Lemma 5.8. Hence we complete the proof. \qed


\end{document}